\documentclass[journal]{IEEEtran}
\ifCLASSINFOpdf
\else
\fi
\newtheorem{theorem}{\textbf{Theorem}}{}
{}
\newtheorem{remark}{\textbf{Remark}}{}
\newtheorem{assumption}{\textbf{Assumption}}{}
{}
\newtheorem{proposition}{\textbf{Proposition}}{}
{}
{}
\usepackage{cite}
\usepackage{amsmath,amssymb,amsfonts}
\usepackage{algorithmic}
\usepackage{graphicx}
\usepackage{textcomp}
\usepackage{xcolor}

\begin{document}
\title{Robust Tracking and Model Following Controller Based on Higher Order Sliding Mode Control and Observation: With an Application to MagLev System}

\author{Siddhartha~Ganguly,
        Manas~Kumar~Bera,
        and~Prasanta~Roy
\thanks{S. Ganguly is with the Department
of Systems $\&$ Control Engineering, Indian Institute of Technology Bomabay, Mumbai, Maharashtra 400076
e-mail: siddhartha@sc.iitb.ac.in.}% <-this % stops a space
\thanks{M.K. Bera is with the Department of Electronics and Instrumentation Engineering, National Institute of Technology Silchar, Silchar, Assam, India 788010 e-mail: manas.bera@gmail.com} \thanks{P. Roy is with the Department
of Electrical Engineering, National Institute of Technology Silchar, Silchar, Assam, India 788010
e-mail: proyiitkgp@gmail.com.}}% <-this % stops a space
\maketitle

\begin{abstract}
This paper deals with the design of robust tracking and model following (RTMF) controller for linear time-invariant (LTI) systems with uncertainties. The controller is based on the second order sliding mode (SOSM) algorithm (super twisting) which is the most effective and popular in the family of higher order sliding modes (HOSM). The use of super twisting algorithm (STA) eliminates the chattering problem encountered in traditional sliding mode control while retaining its robustness properties. The proposed robust tracking controller can guarantee the asymptotic stability of tracking error in the presence of time varying uncertain parameter and exogenous disturbances. Finally, this strategy is implemented on a magnetic levitation system (MagLev) which is inherently unstable and nonlinear. While implementing this proposed RTMF controller for MagLev system, a super twisting observer (STO) is used to estimate the unknown state i.e the velocity of the ball which is not directly available for measurement. It has been observed that the RTMF controller based on STA-STO pair, is not good enough to achieve SOSM for a chosen sliding surface using continuous control. As a remedy, continuous RTMF controller based on STA is implemented with a higher order sliding mode observer (HOSMO). The simulated as well as the experimental results are provided to illustrate the effectiveness of the proposed controller-observers pair.
\end{abstract}


% \fi
\IEEEpeerreviewmaketitle

\section{Introduction}
The main task in the design of a model following controller is to come up with a control algorithm which compels the dynamics of a specified plant to follow the model dynamics. The tracking error between the model and plant outputs should asymptotically converge to zero making sure that the plant outputs follow the model output perfectly. This approach has the complete freedom to specify all the design criteria through model to ensure the minimization of error between model and plant output [1]--[2].\\
In the model based control strategy the mismatch between the model of the plant and the actual plant dynamics exists due to the presence of external disturbances, parametric uncertainties and unmodeled dynamics resulting in declination of performance of the closed loop system. This necessitates the design of various robust controllers. Therefore, the problem of RTMF for a class of dynamical system has gathered good amount of attention amongst many researchers in the recent past [3]--[14].
\par In [3], a nonlinear control scheme is proposed to achieve robust tracking of dynamical signals. Linear state feedback controller is developed for robust tracking and model following of uncertain linear system in [4]. In [5], the problem of RTMF is considered for a class of dynamical systems which contain uncertain nonlinear terms and bounded unknown disturbances. However, these robust state feedback tracking controllers do not produce asymptotic tracking; instead, the so-called practical tracking is achieved.
\par Sliding mode control (SMC) technique has been recognized as one of the efficient tools to design robust controller for systems operating under various uncertainty conditions [1], [15]--[17]. It has many attractive features such as insensitive to matched uncertainties, order reduction, invariance and simplicity in design. The two step design procedure of sliding mode controller consists of design of a sliding manifold which reflects the desired performance during sliding and a discontinuous control law which makes sure that the state trajectories evolving in the sliding manifold [1] stays there for all future time. Based on SMC theory, a class of variable structure tracking controllers has been proposed in [6]--[14] for robust tracking and model following of dynamical signals for a class of uncertain dynamical systems. The main disadvantage of these SMC based strategies are the so called ``chattering effect" which affects the actuator performance in practical implementation [15]--[17]. Moreover, most of these works in [6]--[14] have been tested only in simulation and no hardware implementation has been found so far in any literature.
\par In contrast to that, the controller based on STA is the most suitable choice to eliminate chattering due to continuous nature of control signal [17], [18]. It achieves finite-time convergence of not only the sliding variable but its first derivative as well. In the family of HOSM, the most effective SOSM algorithm is the STA because same algorithm can be used as controller [19], observer known as STO [20], [21] and robust exact differentiator [22]. Controller based on STA generates continuous control signal and consequently compensates chattering which ensures all the properties of first order SMC in addition to that it achieves second order sliding motion (SOSM) in the presence of bounded uncertainties.
\par With this motivation this paper presents a new class of RTMF controllers based on STA and then the proposed control algorithm is implemented in magnetic levitation (MagLev) apparatus (Feedback Instruments Ltd. UK, Model No. 33-210).
\par Industrial application of MagLev continues to grow world wide and successfully implemented for many applications like high-speed MagLev passenger trains, frictionless bearing, superconductor rotor suspension of gyroscopes, rocket-guiding projects and vibration isolation systems in semiconductor manufacturing. The position control of the levitated ball is a challenging control problem due to the fact that the system is inherently nonlinear and open loop unstable. In past many advanced controllers like robust sliding mode control [23]--[28], feedback linearization based control strategy [29]--[31], PID control [32]--[34], adaptive control [35]--[37], generalised proportional integral (GPI) control [38]--[39] have been proposed for compensation of magnetic levitation system.
\par The implementation of the proposed RTMF controller based on STA in MagLev system requires the information of position and velocity of the ball, but direct measurement of ball velocity is not available. For estimation of velocity, the measured position data is numerically differentiated in most of the literature dealing with control law implementation for MagLev system. This process is not suitable from the implementation point of view because differentiation will enhance the high frequency noise in the position measurement. In [29], a nonlinear observer with linear error dynamics is designed to estimate the speed. An model based integral re-constructor is proposed in [39] for online ball velocity estimation. However, these methodologies suffer from a drawback is that the estimators are not robust against the uncertainties. This fact was the motivation to design STO which is insensitive with respect to uncertainties to estimate the ball velocity of MagLev which is complemented with proposed RTMF controller based on STA. But it has been observed that with this RTMF controller based on STA with STO it is not possible to achieve SOSM using continuous control. As a remedy, continuous RTMF controller based on STA is implemented with a HOSMO. With this observation it can be said that the performance of RTMF based on STA with STO as well as with HOSMO has been validated experimentally first time with magnetic levitation apparatus.\\
The contribution can be summarized as follows:
\begin{itemize}
  \item The design of RTMF based on STA is proposed for a class of uncertain LTI systems which will generate continuous control signal, and consequently suppress the effect of chattering.
  \item To validate the performance of this proposed controller, the implementation is done in MagLev system.
  \item To avoid the implementation difficulties, the velocity of the levitated ball is estimated using STO and the RTMF based on STA is designed and implemented based on this estimated velocity. With this method it is not possible to achieve continuous control for chosen sliding surface. As a remedy, continuous RTMF controller based on STA is implemented with a HOSMO.
\end{itemize}
\par The rest of the paper is organized as follows. Section 2 gives some preliminaries and problem statement. Section 3 presents the design steps of RTMF controller based on STA for an uncertain LTI system. In section 4 dynamic model of the MagLev system is derived. In selection 5  the selection of model is discussed along with the design of RTMF controller based on STO and HOSMO for MagLev. In section 6  both simulated and experimental results are presented. Finally, Section 7 concludes the paper.

\section{Problem formulation and assumption}
Consider a $n^{\text{th}}$ order linear time invariant (LTI) uncertain dynamical system as
\begin{align}\label{E1}
\dot{x}(t)&=(A+\Delta A(t,x))x(t)+(B+\Delta B(t,x))u(t)\nonumber \\&+D(t,x)f(t,x)\\
y(t)&=Cx(t)
\end{align}
where $x(t) \in \mathcal{X}\subseteq \mathbb{R}^{n}$ is the state vector, the control input is given by the mapping $\mathbb{R} \ni t \mapsto u(t) \in \mathbb{U} \subset \mathbb{R}^m \neq \varnothing $ which is Lebesgue measureable, $ y(t)\in \mathcal{Y} \subseteq \mathbb{R}^p $ is the output vector which is to track the reference input $y_r(t)\in \mathbb{R}^p$ and $A\in \mathbb{R}^{n\times n}, B\in \mathbb{R}^{n\times m}, C\in \mathbb{R}^{p\times n}$ are constant real matrices. The matrix $D\in\mathbb{R}^{n\times l}$ is known and the function $f(t,x) : \mathbb{R}_{+} \times \mathbb{R}^{n} \rightarrow \mathbb{R}^l$ is unknown exogenous input. Throughout the paper the following assumptions hold.

\begin{assumption}\label{a1}
The pair $(A,B)$ is completely controllable.
\end{assumption}
\begin{assumption}\label{a2}
 The parametric uncertainties viz., $ \Delta A(t,x) $, $ \Delta B(t,x) $, $ D(t,x) $ and the exogenous input signal $ f(t,x) $ are Lipschitz continuous in their arguments. This is indeed necessary to ensure uniqueness of system trajectory.
\end{assumption}

\begin{assumption}\label{a3}
Any uncertainties and disturbances entering into the system satisfy the matching condition, i.e., all uncertain quantities reside in range space of $ B $, $ \mathcal{R}(B) $. Therefore, there exist matrices $ R(t,x): \mathbb{R} \times  \mathbb{R}^{n} \rightarrow \mathbb{R}^{m \times n} $, $ E(t,x): \mathbb{R} \times  \mathbb{R}^{n} \rightarrow \mathbb{R}^{m \times m} $ and $ \mathcal{G}(t,x): \mathbb{R} \times \mathbb{R}^{n} \rightarrow \mathbb{R}^{m \times l} $ such that
\begin{align*}
%\label{E3}
\Delta A(t,x) &= BR(t,x), \quad
\Delta B(t,x) = BE(t,x) \\& {\text {and}}~D(t,x) = B\mathcal{G}(t,x),~ \forall (t,x)\in \mathbb{R}_{+}\times \mathbb{R}^n
\end{align*}

%for every $ (t,x)\in \mathbb{R}_{+} \times \mathbb{R}^{n} $.
\end{assumption}
From the structural assumption all uncertainties can be lumped and the system \eqref{E1} can be rewritten as
\begin{align}\label{E2}
\begin{split}
\dot{x}(t)&=Ax(t)+B\big{(}u(t)+w(t,x)\big{)}\\
y(t)&=Cx(t)
\end{split}
\end{align}
where, $\mathbb{R} \times \mathcal{X} \ni (t,x) \mapsto w(t,x):=R(t,x)x(t)+E(t,x)u(t)+\mathcal{G}(t,x)f(t,x) \in \mathbb{R}^{m}$
is a lumped uncertain vector.
\begin{assumption}\label{a4}
$ w(t,x) \in \mathcal{C}^1\left(\mathcal{X}\subseteq \mathbb{R}^n\right) $ i.e continuously differentiable in $x$ and piecewise continuous in $ t \in \mathbb{R} $. Then there exists a positive constant $\theta_{M} \in \left[0, \infty \right[$ such that
\begin{align}\label{E3}
\|w(t,x)\| \leq \theta_{M}~,~
\forall (t,x)\in \mathbb{R} \times \mathcal{X}\end{align}

\end{assumption}
\par In this paper the objective is to design the control law to make the output of the system \eqref{E2} to follow the output of the reference model. Consider the reference model is given by
\begin{align}\label{E4}
\begin{split}
\dot{x}_{r}(t) &= A_{r}x_{r}(t)\\
 y_{r}(t) &= C_{r}x_{r}(t)
\end{split}
\end{align}
where $x_{r}(t)\in \mathbb{R}^{n_{r}}$ is the state vector of the reference model, $y_{r}(t) \in \mathbb{R}^p$ has the same dimension as $y(t)$ and $A_{r},C_{r}$ are constant matrices of appropriate dimensions. The reference input $y_{r}(t)$ is assumed to be the output of the model. Further, it is assumed that the model states are bounded, i.e., there exists a finite positive constant $ \mathcal{L} $ such that $ \forall t \geq t_{0} $, $ \|x_{r}\| \leq \mathcal{L} $. The reference model represents the ideal response that the controlled system must follow.
\par It has been shown in [4] that if there exist matrices $ G \in \mathbb{R}^{n\times{n_{r}}} $ and $ H \in \mathbb{R}^{m\times{n_{r}}} $ such that the following matrix algebraic relation holds
\begin{align}\label{E7}
\left[\begin{array}{cc}
    A & B \\
    C & 0
  \end{array}\right]
  \left[\begin{array}{c}
           G \\
           H
         \end{array}\right]&=\left[\begin{array}{c}
                    GA_{r} \\
                    C_{r}
                  \end{array}\right]
\end{align}
then the RTMF control based on STA proposed latter in \eqref{control} will ensure that the output of the system \eqref{E2} will follow the output of \eqref{E4}. If a solution to \eqref{E7} is not found, then a different reference model must be chosen. A method of solution to \eqref{E7} is discussed in [4] which is revisited in the following discussion briefly. To solve \eqref{E7}, it is assumed that
\begin{align*}
{\tt rank}\left[\begin{array}{cc}
    A & B \\
    C & 0
  \end{array}\right] &= n+p
\end{align*}
\par This condition is satisfied if the nominal system is controllable and the number of outputs is less than or equal to the number of inputs, i.e., $ p \leq m $. If this satisfy, then solution to \eqref{E7} can be written as
\begin{align}\label{E10}
\left[\begin{array}{c}
           G \\
           H
         \end{array}
       \right] &= \left(\left[
  \begin{array}{cc}
    A & B \\
    C & 0
  \end{array}
\right]^{\top}\left[
  \begin{array}{cc}
    A & B \\
    C & 0
  \end{array}
\right]\right)^{-1} \nonumber \times \\&\left[
  \begin{array}{cc}
    A & B \\
    C & 0 \\
  \end{array}
\right]^{\top}\left[
                  \begin{array}{c}
                    GA_{r} \\
                    C_{r} \\
                  \end{array}
                \right] \nonumber\\
                &=\underbrace{\left[
  \begin{array}{cc}
    \Omega_{11} & \Omega_{12}\\
    \Omega_{21} & \Omega_{22} \\
  \end{array}
\right]}_{=:\Omega}\left[
                  \begin{array}{c}
                    GA_{r} \\
                    C_{r} \\
                  \end{array}
                \right]
\end{align}
%or
%\begin{align}\label{E10}
%\left[
%         \begin{array}{c}
%           G \\
%           H \\
%         \end{array}
%       \right]=\left[
%  \begin{array}{cc}
%    \Omega_{11} & \Omega_{12}\\
%    \Omega_{21} & \Omega_{22} \\
%  \end{array}
%\right]\left[
%                  \begin{array}{c}
%                    GA_{m} \\
%                    C_{m} \\
%                  \end{array}
%                \right]
%\end{align}
where $ \Omega_{11} $ is an $ n\times n $ matrix, $ \Omega_{12} $ is an $ n\times p $ matrix, $ \Omega_{21} $ is an $ m\times n $ matrix, and $ \Omega_{22} $ is an $ m\times p $ matrix. Specifically, if $ p = m $ then :
\begin{align}\label{E11}
\Omega:=\left[
  \begin{array}{cc}
    A & B \\
    C & 0 \\
  \end{array}
\right]^{-1}\nonumber
\end{align}
and $ \Omega $ is partitioned accordingly in the same way. An equivalent form of \eqref{E10} is
%\begin{align}
%    G &= \Omega_{11}GA_{m} + \Omega_{12}C_{m} \\
%    H &= \Omega_{21}GA_{m} + \Omega_{22}C_{m}
%\end{align}
\begin{equation}\label{E11.1}
  G = \Omega_{11}GA_{r} + \Omega_{12}C_{r}
\end{equation}
\begin{equation}\label{E11.2}
  H = \Omega_{21}GA_{r} + \Omega_{22}C_{r}
\end{equation}
\par The above relations can be solved simultaneously for $ G $ and $ H $. Though there are many methods exist to solve for the above relation here we discuss the method given in [4]. Let $ M $ be any $ n\times m $ matrix
\begin{align*}
M &:= \begin{bmatrix}
           m_{1} \\
           m_{2} \\
          \vdots\\
          m_{n}
       \end{bmatrix}
\end{align*}
where $m_{i}$ denotes the $i^{ \text {th}}$ row of $ M $. We define a stacking operator $ \Phi(M) $ for the matrix $ M $ as
\begin{align}\label{E13}
\Phi(M) & := \begin{bmatrix}
           m^{\top}_{1}  \\
           m^{\top}_{2} \\
           \vdots \\
           m^{\top}_{n}
       \end{bmatrix}
\end{align}
The vector is obtained by stacking the transpose of the rows of $ M $. To solve \eqref{E11.1}, we write it as
\begin{equation}\label{E14}
G-\Omega_{11}GA_{r}=\Omega_{12}C_{r}
\end{equation}
This equation can be solved using the Kronecker product $\otimes$. The Kronecker product of two matrices $ U \in \mathbb{R}^{n \times m} $ and $ V \in \mathbb{R}^{q \times r}$ is defined as
\begin{align}\label{E15}
U \otimes V &:= \begin{bmatrix}
        u_{11}V & u_{22}V & \dotsb & u_{1m}V \\
        \vdots &  \vdots & \ddots & \vdots\\
         u_{n1}V & u_{n2}V & \dotsb & u_{nm}V
       \end{bmatrix}
\end{align}
The resultant matrix $ U \otimes V $ is of the order $ nq \times mr $. Using the stacking operator and the Kronecker product, \eqref{E14} reduces to
\begin{align}\label{E16}
[I_{n}\otimes I_{n_{r}} - \Omega_{11} \otimes A^\top_{r}]\Phi(G)=\Phi(\Omega_{12}C_{r})
\end{align}
which is a system of linear equations for $ nm $ components of $ \Phi(G) $. Then $ G $ is found by unstacking $ \Phi(G) $. Once $ G $ is solved, substituting for $ G $ into \eqref{E11.2} yields $ H $.

%\IEEEpubidadjcol

\section{Design of robust tracking controller}
Our design objective can be stated as: given a desired model and its output, find a RTMF controller to force the uncertain system under consideration to behave as the ideal model, i.e., have the same or similar rise time, settling time, damping, etc. To achieve this, we propose a RTMF controller based on STA which can guarantee the output $ y(t) $ of uncertain system \eqref{E2} will follow the output $ y_{r}(t) $ of reference model \eqref{E4} ensuring asymptotic convergence of tracking error to zero. Let the tracking error be defined as
\begin{align}
e(t) &:= y(t)-y_{r}(t)
\end{align}
Then the tracking control law is proposed as
\begin{align}\label{control}
u(t) &:= Hx_{r}(t)+v(t)
\end{align}
where $ v(t) \in \mathbb{U}\subseteq \mathbb{R}^m $ is a super twisting control law. A new auxiliary state vector $ z(t) \in \mathcal{X}$ is defined as follows
\begin{align}\label{z}
\mathbb{R} \ni t \mapsto z(t):=x(t)-Gx_r(t) \in \mathcal{X}\subseteq \mathbb{R}^n
\end{align}
\par Using \eqref{E7} and \eqref{z}, the relation between tracking error $e(t)$ and new state vector $z(t)$ can be written as
\begin{align}\label{E17}
e(t) &= Cz(t)
\end{align}
From \eqref{E17}, we can see that $ \|e(t)\| \leq \|C\| \|z(t)\| $, and as $\|C\|<\infty$, it's evident that $ \|z(t)\| \to 0 \implies $   $ \|e(t)\| \to 0$.
%\begin{align}\label{E18}
%\|z(t)\|\rightarrow 0~~\text{implies}~~\|e(t)\|\rightarrow 0
%\end{align}
So, it's enough to only think about stability of $z(t)$.
From \eqref{E4}, \eqref{control} and \eqref{z}, the auxiliary system can be expressed as
\begin{align}\label{system1}
\dot{z}(t)&=Az(t)+Bv(t)+Bw(t, z(t),x_r(t))
\end{align}
\par A convenient way to solve this problem is to first transform the system \eqref{system1} into a suitable canonical form. In order that, partition \eqref{system1} as
\begin{align*}
B=\left[
    \begin{array}{c}
      B_{1 }\\
      B_{2} \\
    \end{array}
  \right]
\end{align*}
where $ B_{1} \in \mathbb{R}^{(n-m)\times{m}} $ and $ B_{2} \in \mathbb{R}^{m\times{m}} $, and $ \det B_{2}\neq 0 $.
Define the following orthogonal transformation
\begin{align}\label{transform}
T_{1} &:= \left[\begin{array}{cc}
    \mathbb{I}_{n-m} & -B_{1}B_{2}^{-1} \\
    0 & B_{2}^{-1 }
  \end{array}
  \right]
\end{align}
where $\mathbb{I}_{n-m}$ is a $(n-m) \times (n-m)$ indentity matrix, now define the new coordinate of transformation as
\begin{align}\label{tran}
\left[\begin{array}{c}
    \eta \\
    \xi \\
  \end{array}\right] &:= T_{1}z(t)
\end{align}

Therefore, \eqref{system1} can be written with the help of above transformation in the regular form as
 \begin{align}\label{system3}
%\begin{split}
\left[
  \begin{array}{c}
    \dot{\eta} \\
    \dot{\xi} \\
  \end{array}
\right]=\left[
          \begin{array}{cc}
            A_{11} & A_{12} \\
            A_{21} & A_{22} \\
          \end{array}
        \right]\left[
                 \begin{array}{c}
                   \eta \\
                   \xi \\
                 \end{array}
               \right]+\left[
                         \begin{array}{c}
                           0 \\
                           \mathbb{I}_{m}\\
                         \end{array}
                       \right]v(t)+\left[
                         \begin{array}{c}
                           0 \\
                           \mathbb{I}_{m}\\
                         \end{array}
                       \right]w(t)
%\end{split}
\end{align}
where $\eta \in \mathbb{R}^{n-m}$  and $\xi \in \mathbb{R}^{m}$. In order to achieve robustness against matched uncertainty, we design STA for the \eqref{system3}. Since STA is a first order control algorithm, so a hypersurface is designed such that the control appears in its first derivative. This technique achieves asymptotic stability of the system once this hypersurface is reached in finite time.
\par We define the sliding variable as 
\begin{align}\label{E19}
\sigma:=\left[
    \begin{array}{cc}
      -K & \mathbb{I}_{m} \\
    \end{array}
  \right]\left[
            \begin{array}{c}
              \eta \\
              \xi \\
            \end{array}
          \right]
\end{align}
where $ K \in \mathbb{R}^{m \times (n-m)} $ is chosen by some appropriate design procedure and $ \sigma^{\top} = \begin{bmatrix}\sigma_{1} & \dotsb & \sigma_{m}\end{bmatrix} $. Further, define a linear change of coordinates by
\begin{align*}
%\label{E20}
T_{2}:= \left[
                  \begin{array}{cc}
                    \mathbb{I} & 0 \\
                    -K & \mathbb{I} \\
                  \end{array}
                \right]
\end{align*}
If this transformation $ T_2 $ represents the new coordinates as
\begin{align*}
%\label{E21}
\left[
  \begin{array}{c}
    \eta \\
    \sigma \\
  \end{array}
\right]:=T_{2}\left[
  \begin{array}{c}
    \eta \\
    \xi \\
  \end{array}
\right]
\end{align*}
then dynamics in the new coordinates is
\begin{align}\label{E22}
\left[
  \begin{array}{c}
    \dot{\eta }\\
   \dot{ \sigma }\\
  \end{array}
\right]=\mathbb{A}
\left[
  \begin{array}{c}
    \eta \\
    \sigma \\
  \end{array}
\right]+\left[
                         \begin{array}{c}
                           0 \\
                           \mathbb{I}_{m}\\
                         \end{array}
                       \right]v(t)+\left[
                         \begin{array}{c}
                           0 \\
                           \mathbb{I}_{m}\\
                         \end{array}
                       \right]w(t)
\end{align}
where
\begin{align*}
\mathbb{A}:=\left[
          \begin{array}{cc}
            A_{11}+A_{12}K & A_{12} \\
            \underbrace{A_{21}+A_{22}K-K(A_{11}+A_{12}K)}_{=:\mathcal{A}} & A_{22}-KA_{12} \\
          \end{array}
        \right]
\end{align*}
Design the control input $ v(t) $ as
\begin{align}\label{E231}
v(t) &= -\left[
          \begin{array}{cc}
           0 & 0 \\
            \mathcal{A} & A_{22}-KA_{12} \\
          \end{array}
        \right]\left[
  \begin{array}{c}
    \eta \\
    \sigma \\
  \end{array}
\right]+v'(t)
\end{align}
Substituting the control input in \eqref{E22} the system equation becomes
\begin{align}\label{E23}
\left[
  \begin{array}{c}
    \dot{\eta }\\
   \dot{ \sigma }\\
  \end{array}
\right] &= \left[
          \begin{array}{cc}
            A_{11}+A_{12}K & A_{12} \\
            0 & 0 \\
          \end{array}
        \right]
\left[
  \begin{array}{c}
    \eta \\
    \sigma \\
  \end{array}
\right]+\left[
                         \begin{array}{c}
                           0 \\
                           \mathbb{I}_{m}\\
                         \end{array}
                       \right]v'(t)\nonumber\\
&\quad +\left[
                         \begin{array}{c}
                           0 \\
                           \mathbb{I}_{m}\\
                         \end{array}
                       \right]w(t)
\end{align}
where the components of $ v'(t) $ are
\begin{align}\label{E24}
\begin{split}
v_{i}'(t)&=-k_{1}|\sigma_{i}|^{\frac{1}{2}}\text{sign}(\sigma_{i})+\Omega_{i}\\
\dot{\Omega}_{i}&=-k_{2}\text{sign}(\sigma_{i})
\end{split}
\end{align}
where $ k_{1}, k_{2} $ are the scalar gains. Let $ L_{i} := \Omega_{i}+w_{i} $ and substituting \eqref{E24} in \eqref{E23} then the closed loop system
\begin{align}
%\begin{split}
\dot{\eta}&=(A_{11}+A_{12}K)\eta+A_{12}\sigma \label{E24a}\\
%\dot{\sigma_{i}}&=-k_{1}\text{sign}(\sigma_{i})|\sigma_{i}|^{\frac{1}{2}}+L_{i}\label{E24b}\\
%\dot{L_{i}}&=-k_{2}\text{sign}(\sigma_{i})+\dot{w}_{i}(t)\label{E24c}
\dot{\sigma}&=-k_{1}\lfloor\sigma\rceil^{\frac{1}{2}}+L\label{E24b}\\
\dot{L}&=-k_{2}\text{Sign}(\sigma)+\dot{w}(t)\label{E24c}
%\end{split}
\end{align}
where the vector $ \lfloor\sigma\rceil^{\frac{1}{2}} = \begin{bmatrix}\lfloor\sigma_1\rceil^{\frac{1}{2}} & \dotsb & \lfloor\sigma_m\rceil^{\frac{1}{2}}\end{bmatrix}^\top $ with each component $ \lfloor\sigma_i\rceil^{\frac{1}{2}} = |\sigma_{i}|^{\frac{1}{2}}\text{sign}(\sigma_{i}) $ for $ i = 1, \dotsc, m $ and $ \text{Sign}(\sigma) = \begin{bmatrix}\text{sign}(\sigma_1) & \dotsb & \text{sign}(\sigma_m)\end{bmatrix}^\top $. The closed loop dynamics \eqref{E24a}--\eqref{E24c} consist of a set of differential inclusion and the solution to this are understood in sense of Filippov [40].
\begin{theorem}
Consider the system \eqref{E24a}-\eqref{E24c}. Then the closed loop dynamics \eqref{E24b}-\eqref{E24c} is finite time stable if the gains are selected such that $k_{1}>0$ and $k_{2}>\max_{i}|\dot{w}_{i}(t)|$. Moreover, the reduced dynamics \eqref{E24a} is asymptotically stable.
\end{theorem}
%if and only if there exists gain $K$ such that $A_{11}+A_{12}K$ is Hurwitz.

% \begin{proof}

\textbf{Proof : }For any $ i^{\text {th}} $ input, the dynamics \eqref{E24b}-\eqref{E24c} can be rewritten as
\begin{align}
%\begin{split}
\dot{\sigma_{i}}&=-k_{1i}|\sigma_{i}|^{\frac{1}{2}}\text{sign}(\sigma_{i})+L_{i}\label{E25b}\\
\dot{L_{i}}&=-k_{2i}\text{sign}(\sigma_{i})+\dot{w}_{i}(t)\label{E25c}.
%\end{split}
\end{align}
%Then there exists a time $t_{1}$ such that for all $t \geq t_{1}$
%\begin{align}
%\sigma_{1}(t)=\sigma_{2}(t)=\dotsb=\sigma_{m}(t)=0.
%\end{align}
\par To show \eqref{E25b} and \eqref{E25c} finite time stable, the Lyapunov function can be chosen as $ V_{i}(\Psi) = \Psi_{i}^{\top}P_{i}\Psi_{i} $ where $ \Psi_{i}^{\top} = \left[|\sigma_{i}|^{\frac{1}{2}}\text{sign}(\sigma_{i})~~L_{i}\right]$,  then the trajectories of the system \eqref{E25b}-\eqref{E25c} will converge to origin in a finite time smaller than $ t_{i} $ [41]
\begin{align}
t_{i}=\frac{2}{\gamma_{i}}V_{i}^{\frac{1}{2}}(\Psi_{i}(0))\nonumber
\end{align}
where
\begin{align}
\gamma_{i}=\frac{\lambda_{\min}^{1/2}\{P_{i}\}\lambda_{\min}\{Q_{i}\}}{\lambda_{\max}\{P_{i}\}}\nonumber
\end{align}
for positive and symmetric definite matrices $P_{i}$ and $Q_{i}$. Here the gains $k_{ji}$ for $j=1,2$ can be selected to achieve $\sigma_{i}$ and $L_{i}$ zero in some finite time.
\par Select the gains $k_{1}$ and $k_{2}$ as $ k_1 = \min_{i}|k_{1i}| $ and $ k_2 = \max_{i}|k_{2i}| $. The gains $ k_1 $ and $ k_2 $ are enough to bring $ \sigma \equiv 0 $ in finite time.
%\begin{align*}
%k_1 = \min_{i}|k_{1i}|, \quad &k_2 = \max_{i}|k_{2i}|.
%\end{align*}
Hence from \eqref{E24a}, the reduced order dynamics or zero dynamics of the system with respect to $ \sigma $ becomes
\begin{align}\label{reduced}
\dot{\eta}&=(A_{11}+A_{12}K)\eta.
\end{align}
The pair $ (A_{11},A_{12}) $ is controllable pair and $ K $ can be designed in such way that $\text{spec}\left(A_{11}+A_{12}K \right) \in \mathbb{C}^-$ i.e \eqref{reduced} becomes globally asymptotically stable. This ends the proof. $\blacksquare$
% \end{proof}
\begin{proposition}
Consider the model following problem of uncertain system $\eqref{E1}$ satisfying Assumptions \ref{a1}--\ref{a4}. Then, the control law $ u(t) $ given by \eqref{control}, \eqref{E231} and \eqref{E24} guarantees asymptotically stability of the tracking error $ e(t) $.
\end{proposition}
% \begin{proof}
\textbf{Proof : }From Theorem 1, it is shown that the sliding surface $ \sigma^{\top} = [\sigma_{1} \dotsb \sigma_{m}] $ goes to zero in finite time. Using \eqref{reduced}, it can be concluded that
\begin{align}\label{E25}
\begin{split}
\lim_{t\to\infty} \eta(t) &= 0
\end{split}
\end{align}
Now using \eqref{E19} and \eqref{E25}, it yields
\begin{align}\label{E26}
\lim_{t\to\infty} \xi(t) &= 0
\end{align}
Then, it follows immediately from \eqref{tran}, \eqref{E25} and \eqref{E26} that
\begin{align}
\lim_{t\to\infty} z(t) &= 0
\end{align}
Recalling relationship between $ e(t) $ and $ z(t) $, i.e. $ e(t) = Cz(t) $, we obtain that the tracking error $ e(t) $ also decreases asymptotically to zero. $\blacksquare$
% \end{proof}
\par The real time implementation of the proposed control law for magnetic levitation system is developed in the following sections.
\section{Dynamical Model of the MagLev system}\label{sec4}
\begin{figure} [h]
\centering
\includegraphics[width=3in]{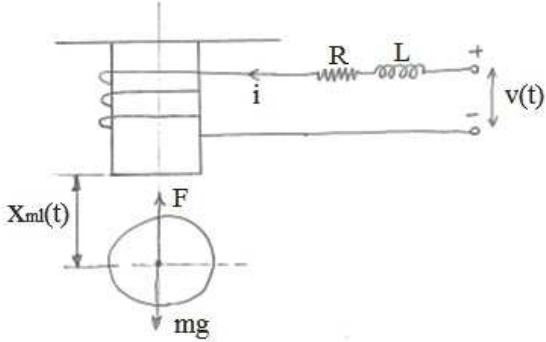}
\caption{Schematic of magnetic levitation system }
\label{MAG}
\end{figure}

A simplified schematic diagram of the MagLev system is shown in Fig.\ref{MAG}. The system is mainly computer controlled and user's specific control law can be implemented and applied on the system in realtime using MATLAB/SIMULINK environment
through a  PCI 1711 AD card.\\

% \begin{figure} [h]
% \centering
% \includegraphics[width=2.5in]{maglel_setup.eps}
% \caption{The experimental setup of the MagLev system }
% \label{MAG_SETUP}
% \end{figure}
 Our main task here is vertical position control of the levitated ball by application of a suitable current i(t) by means of applying a voltage v(t) in the coil. R and L represents the resistance and inductance of the coil, respectively. Here $x_{\text{ml}}(t)$ denotes the distance by which the coil magnet and the ball’s center are separated. The  nonlinear force applied by the magnet is $\mathcal{K} \frac{i^2}{x^2_{\text{ml}}}$, where $\mathcal{K}$ is a positive constant. The differential equation capturing dynamics of the system is of the form
\begin{align}
m\ddot{x} &= mg -\mathcal{K} \frac{i^2}{x^2_{\text{ml}}} \label{E31}\\
i &= {\mathcal{K}}_1u \label{E32}
\end{align}
where $m$ denotes mass of the levitated ferromagnetic ball, $g$ is the gravitational acceleration and ${\mathcal{K}}_1$ is a  proportionality constant [42] which depends upon the input voltage and coil current. To design the proposed robust model following controller, the model is linearized around an equilibrium point $(x_0,i_0)$ and the transfer function of the system is found as
\begin{align}\label{E33}
\frac{ x_{\text{ml}}(s)}{i(s)} &=\frac{-{\mathcal{K}}_i}{s^2+{\mathcal{K}}_x}
\end{align}
where ${\mathcal{K}}_i=\frac{2mg}{i_{0}}$, $\mathcal{K}_{x}=-\frac{2mg}{x_{0}}$. The ball position  $x_{\text{ml}}$ is measured by the IR sensor in voltage $x_{v}$ and they are linearly related with a proportionality constant $\mathcal{K}_{2}$. Finally with $x_{v}$, the above transfer function can be rewritten as
\begin{align}\label{E333}
\frac{ x_{v}(s)}{u(s)} =\frac{-\mathcal{K}_{1}\mathcal{K}_{2}{\mathcal{K}}_i}{s^2+{\mathcal{K}}_x}=\frac{-3518.85}{s^{2}-2180}
\end{align}
The nominal system parameters are shown in Table \ref{tab1}. From the open loop transfer function it is clear that the system is severely unstable.
\par In state space representation \eqref{E333} can be written as
\begin{align} \label{E35}
	 \begin{split}
        \dot{x}(t)&=\left[
                      \begin{array}{cc}
                        0 & 1 \\
                        2180 & 0 \\
                      \end{array}
                    \right]x(t)+\left[
                              \begin{array}{cc}
                                0  \\
                                -3518.85 \\
                              \end{array}
                            \right]u(t)\\
y(t)&=\left[
                              \begin{array}{cc}
                                1 & 0
                              \end{array}
                            \right]x(t).
      \end{split}
 \end{align}
From \eqref{E35}, it is clear that the MagLev is a second order electro-mechanical system whose relative degree is two [43].
% \begin{table}[!b]
%\processtable{System parameters\label{tab1}}
%{\begin{tabular*}{20pc}{@{\extracolsep{\fill}}ll@{}}\toprule
%Parameter & Value \\

%\midrule
%$ m $ & $ 0.02 $ Kg   \\
%$ g $ & $ 9.81 $  $ \text{m}$/$\text{s}^2$   \\
%$ i_0 $ & $ 0.8 $ A   \\
%$ x_0 $ & $ 0.009 $ m \\
%$ {k}_1 $ & $ 1.05 $ A/V  \\
%$ {k}_2 $ & $ 143.48 $ V/m \\
%$ u $ & $ \pm 5 $ V \\
%$ x_{v} $ & $ -1.5 $ V $\text{to}$ $ 3.75 $ V \\

%\botrule
%\end{tabular*}}{}
%\end{table}

\begin{table}[!t]
%\renewcommand{\arraystretch}{1.3}
 %if using array.sty, it might be a good idea to tweak the value of
% \extrarowheight as needed to properly center the text within the cells
\caption{System parameters}
\label{tab1}
\centering
% Some packages, such as MDW tools, offer better commands for making tables
% than the plain LaTeX2e tabular which is used here.
\begin{tabular}{|c||c|}
\hline
Parameter & Value \\
\hline
$ m $ & $ 0.02 $ Kg \\
\hline
$ g $ & $ 9.81 $  $ \text{m}$/$\text{s}^2$   \\
\hline
$ i_0 $ & $ 0.8 $ A   \\
\hline
$ x_0 $ & $ 0.009 $ m \\
\hline
$ {\mathcal{K}}_1 $ & $ 1.05 $ A/V  \\
\hline
$ {\mathcal{K}}_2 $ & $ 143.48 $ V/m \\
\hline
$ u $ & $ \pm 5 $ V \\
\hline
$ x_{v} $ & $ -1.5 $ V $\text{to}$ $ 3.75 $ V \\
\hline
\end{tabular}
\end{table}
 \section{Design of the proposed controller for MagLev system}
In this section the proposed robust tracking and model following control strategy is designed for the MagLev system which is affected by model uncertainties and external disturbances, then \eqref{E35} can be rewritten as
\begin{align} \label{E36}
	 %\begin{split}
        \dot{x}(t)&=\left[
                      \begin{array}{cc}
                        0 & 1 \\
                        2180 & 0 \\
                      \end{array}
                    \right]x(t)+\left[
                              \begin{array}{cc}
                                0  \\
                                -3518.85 \\
                              \end{array}
                            \right]\big{(} u(t)+ w(t) \big{)} \\\label{E361}
y(t)&=\left[
                              \begin{array}{cc}
                                1 & 0
                              \end{array}
                            \right]x(t)
      %\end{split}
 \end{align}
where the lumped disturbance applied to the plant is chosen as $w(t):=5 \sin(t)$.
\par The task of the robust model following controller is to ensure the asymptotic convergence of tracking error between model and the plant output to zero. The essential part for the controller design is to select a suitable model which the plant needs to follow. The selection of a suitable model for the MagLev system is  discussed in the following subsection.

\subsection{Selection of the Model}
\begin{figure}[h!]
  \centering
 \includegraphics[width=3.8in]{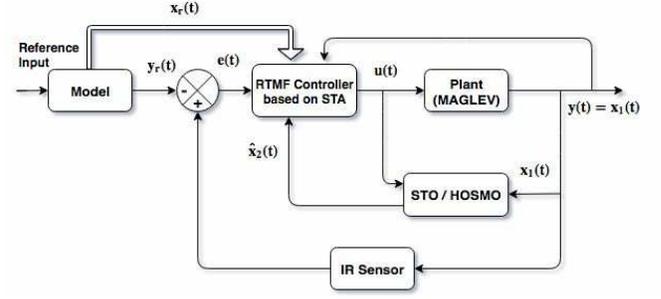}
 \caption{Block diagram}
\label{fig1}
\end{figure}

Proportional-integral-derivative (PID) controller is widely used for the effective regulation of the real systems without uncertainty. In this work, the nominal MagLev system with PID controller is used to find a suitable model. With PID controller $G_{c}(s):=K_{p}+\frac{K_{i}}{s}+K_{D}s$ and the open loop transfer function of MagLev $G_{p}(s):=\frac{-C_{1}}{s^{2}-C_{2}}$, where $C_{1}=3518.85$, $C_{2}=2180$, the closed loop transfer function is computed as follows
\begin{align}\label{E33}
T_{CL}(s)&:=\frac{G_{c}(s)G_{p}(s)}{1+G_{c}(s)G_{p}(s)}\nonumber\\
&=\frac{-(K_{d}s^{2}+K_{p}s+K_{i})C_{1}}{s^{3}-{K}_{d}C_{1}s-(K_{p}C_{1}+C_{2})s-K_{i}C_{1}}
\end{align}\\
The values of the controller gains $ \left\{ K_p=-4.8 , K_i =-98, K_d= -0.06\right\}$ are selected such that the closed loop poles are located at $-70$ [39]. With these gains, \eqref{E33} can be represented as follows
\begin{align}
\dot{x}_{r}(t) &= A_{r}x_{r}(t) \label{E38} \\
y_{r}(t)&=C_{r}x_{r}(t) \label{E39}
\end{align}
where
$A_{r}:=\left[ \begin{array}{ccc} 0 & 1 & 0 \\ 0 & 0 & 1 \\ -343000 & -14700 &  -210  \end{array} \right]$  and \\\\ $C_{r}:=\left[\begin{array}{ccc}343000 & 0 & 0\end{array}\right] $ \\
Now from \eqref{E7}, $G$ and $H$ matrices can be easily computed with the knowledge of the system and model matrices which is shown below
\begin{align} \label{G}
	 %\begin{split}
        G&:=\left[
            \begin{array}{ccc}
              343000 & 0 & 0 \\
              0 & 343000 & 0 \\
            \end{array}
          \right]\\
H&:=\left[ \begin{array}{ccc}212500 & 0 & -100 \end{array}\right] \label{H}
      %\end{split}
 \end{align}
To apply the proposed RTMF controller based on STA on MagLev which is a second order electromechanical system with relative degree two to compensate the chattering problem, the design of sliding surface \eqref{E19} will be such that it's relative degree must be one. In that situation we require the information of two states i.e. ball position and ball velocity of MagLev system. In practice only ball position is available for measurement. In the following section it has been shown that using STO it is possible to estimate the ball velocity in finite time in the presence of uncertainties and using this estimated information, the proposed RTMF controller based on STA is possible to design. The block diagram of this proposed scheme is shown in Fig.~\ref{fig1}.

%%%%%%%%%%%%%%%%%%%%%%%%%%%%%%%%%%%%%%%%%%%%%%%%%%%%%%%%%%%%%%%%%%%%%%%%%%%%%%%%%%%%%%%%%%%%%%%%%%%%%%%%%%%%%%%%%%

\subsection{Design of Proposed RTMF Controller based on STO}
To design the STO to estimate the states of MagLev system \eqref{E36}, let us rewrite the plant dynamics \eqref{E36}
\begin{align}
\dot{x}_1(t)&=x_2(t)\\ \nonumber
\dot{x}_2(t)&=2180x_1(t)-3518.85u(t)+w(x,t)\\ \nonumber
y(t)&=x_1(t)
\end{align}
The STO dynamics can be written as
\begin{align}
\dot{ \hat{x}}_1 &=\hat{x}_2+K_1 \text{sign}(e_1) \label{E46} \\
\dot{\hat{x}}_2 &=K_2\text{sign}(e_1)-3518.85\,u(t)+2180\hat{x}_1 \label{E47}
\end{align}
Where $e_1:=x_1-\hat{x}_1$ and $e_2:=x_2-\hat{x}_2$ are the estimation errors. Then the error dynamics can be written as
\begin{align}
    \dot{e}_1&=e_2-K_1|e_1|^{1/2}\text{sign}(e_1) \label{E48} \\
    \dot{e}_2&=-K_2\text{sign}(e_1)+2180e_1+ w(t,x) \label{E49}
\end{align}
From Assumption 4, $\|w(t,x)\| \leq \theta_{M}$. In [41], [44] the finite time stability of the above error dynamics has been established  with the choice of $K_1=1.5\sqrt{\theta_{M}}$ and $K_{2}=1.1\theta_{M}$. This ensures the finite time convergence of estimation errors $e_{1}$ and $e_{2}$ to zero which essentially leads to $x_1=\hat{x}_1$ and $x_2=\hat{x}_2$ in finite time.
\par To design the proposed RTMF controller, system is required to transform in $z$ coordinate using \eqref{z}. Using the relation $z=x-Gx_r$ for MagLev system following can be written
\begin{equation}
\left[\begin{array}{cc}z_1\\
\hat{z}_2\\\end{array}\right]=\left[\begin{array}{cc}x_1\\
\hat{x}_2\\\end{array}\right]-G\left[\begin{array}{ccc}x_{{r}_1}\\
x_{r_{2}} \\ x_{r_3} \end{array}\right] \label{E50}
\end{equation}
Where $x_{r_1},x_{r_2}$ and $x_{r_3}$ are the states of the chosen model \eqref{E38} and $\hat{x}_2$ is the estimated ball velocity by STO.
Using \eqref{E36}, \eqref{G}, \eqref{E47} and \eqref{E50} the derivative of ${z}_1$ and $\hat{z}_2$ can be written as
\begin{align}
   \dot{z}_1 &=\dot{x}_1-343000\dot{x}_{r_1}=x_2-343000\,x_{r_2} \label{E51}\\
   \dot{ \hat{z}}_2 &=\dot{\hat{x}}_2-343000\dot{x}_{r_2} \nonumber\\
   &=K_2\text{sign}(e_1)-3518.85u(t)+2180\hat{x}_1-343000x_{r_3}\label{E52}
\end{align}
Now choosing a sliding manifold of the following form
\begin{equation}
 \hat{s}=c_1z_1+\hat{z}_2;~~\text{where}~c_{1}>0 \label{E53}
\end{equation}
and differentiating \eqref{E53}, the sliding dynamics can be written as
 \begin{equation}
     \dot{\hat{s}}=c_1\dot{z}_1+\dot{\hat{z}}_2 \label{E54}
 \end{equation}
Now using \eqref{E51} and \eqref{E52}, we can write \eqref{E54} as
 \begin{align}
     \dot{\hat{s}}&=c_1x_2-343000\,c_1x_{r_2}+K_2\text{sign}(e_1) \nonumber \\&-3518.85\,u(t)+2180\, \hat{x}_1-343000x_{r_3} \label{E55}
 \end{align}
Using \eqref{H}, the proposed control law \eqref{control} can be written as
\begin{align}
u(t) &=\left[ \begin{array}{ccc}212500 & 0 & -100 \end{array}\right]\left[\begin{array}{ccc}x_{r_1}\\ x_{r_2}\\ x_{r_3} \\\end{array}\right]+v(t) \label{E56}
\end{align}
Now using \eqref{E56}, replacing $ x_2=\hat{x}_2+e_2 $, and choosing $c_1=1 $ \, we can write \eqref{E55} as
\begin{align}
\dot{\hat{s}} &=\left(\hat{x}_2+e_2 \right) -343000x_{r_2}-7.47 \times 10^8 \,x_{r_1} \nonumber \\& +8885\,x_{r_3}+K_2\text{sign}(e_1)-3518.85\,v(t) +2180\, \hat{x}_1 \label{E57}
\end{align}
Now choosing the control input $v(t)$ as
\begin{align}
v(t)&= -\frac{1}{3518.85} \bigg{(} - \hat{x}_2+343000  x_{r_2} +7.47 \times 10^8 x_{r_1} \nonumber \\& -8885x_{r_3} -K_2\text{sign}(e_1)-\lambda_1 |\hat{s}|^{1/2}\text{sign}(\hat{s}) \nonumber \\& -\int_0^t \lambda_2 \text{sign}(\hat{s}) \mathrm{d}\tau  \bigg{)} \label{E59}
\end{align}
where $\lambda_1>0$ and $\lambda_2>0$ are controller design parameters, \eqref{E57} can be written as
\begin{equation}
\dot{\hat{s}}=e_2-\lambda_1|\hat{s}|^{1/2}\text{sign}(\hat{s})-\int_0^t \lambda_2 \text{sign}(\hat{s})\,\mathrm{d}\tau \label{E59.1}
\end{equation}
Using the fact $ \dot{z}_1=\hat{z}_2 $ , \eqref{E53} and \eqref{E59.1}, the closed loop system can be represented in ($z_1$ , $\hat{s}$) coordinate as follows
 \begin{align*}
 \dot{z}_1&=-z_1+\hat{s}\\
  \dot{\hat{s}}&=e_2-\lambda_1|\hat{s}|^{1/2}\text{sign}(\hat{s})-\int_0^t \lambda_2 \text{sign}(\hat{s})\,\mathrm{d}\tau 
 \end{align*}
So the overall closed loop system with controller-observer together now can be represented as
\begin{align}
\Phi :& \begin{cases}
    \dot{z}_1=-{z}_1+\hat{s}  & \\
    \dot{\hat{s}}=e_2-\lambda_1|\hat{s}|^{1/2}\text{sign}(\hat{s})-\int_0^t \lambda_2 \text{sign}(\hat{s})\,\mathrm{d}\tau &
\end{cases} \nonumber \\
\Lambda :& \begin{cases}
   \dot{e}_1=e_2-K_1|e_1|^{1/2}\text{sign}(e_1) & \\
    \dot{e}_2=-K_2\text{sign}(e_1)+2180e_1+w(t,x) &
\end{cases} \label{E61}
\end{align}
The observer error for the system $\Lambda$ goes to zero in finite time. In [45] it has been established that, the trajectories of the system $\Phi$ cannot escape in finite time. Generally, observer gains are chosen in such a way that the estimation error converges faster. This fact allows to design
the observer and the control law separately, i.e., the separation principle is satisfied. With $e_{2}=0$, the closed-loop system can be written as
\begin{align}\label{E100}
\begin{split}
    \dot{z}_1&=\hat{s}-z_1  \\
    \dot{\hat{s}}&=-\lambda_1|\hat{s}|^{1/2}\text{sign}(\hat{s})+\mu \\
    \dot{\mu}&=-\lambda_2\text{sign}(\hat{s})
\end{split}
\end{align}
From the last two equations of \eqref{E100} it is evident that it has the same structure as STA. So it is easy to conclude that in finite time $\hat{s}=\dot{\hat{s}}=0$, then the reduced order system dynamics can be represented as
\begin{gather*}
\dot{z}_1=-z_1 \\
\hat{z}_2=-{z}_1
\end{gather*}
Therefore both states $z_1$ and $\hat{z}_2$ are asymptotically stable and using Proposition 1 it can be said that the tracking error $e(t)$ also converges to zero asymptotically.
\begin{remark}
The additional term $K_2\text{sign}(e_1)$ is added in $v(t)$ \eqref{E59} because in [21] it has been shown that it is not possible to achieve the continious control when super twisting controller is implemented with STO. To solve this implementation issues, a continuous RTMF controller based on STA is proposed with higher order sliding mode observer (HOSMO) like in [21] that achieves the second order sliding mode.
\end{remark}
\subsection{\emph{Design of Proposed RTMF Controller based on HOSMO}}
To estimate the velocity of the uncertain MagLev system \eqref{E36}, the dynamics of HOSMO [21] is presented as
\begin{align}
\dot{\hat{x}}_1 &= \hat{x}_2 + L_1 |e_1|^{\frac{2}{3}} \text{sign}(e_1) \nonumber \\
\dot{\hat{x}}_2 &= L_2 |e_1|^{\frac{1}{3}} \text{sign}(e_1)+2180 \hat{x}_1-3518.85\,u(t)+\hat{x}_3 \nonumber \\
\dot{\hat{x}}_3 &= L_3\text{sign}(e_1) \label{E67}
\end{align}
where $e_1:=x_1 -\hat{x}_1$ and $e_2:=x_2-\hat{x}_2$ are the estimation error variables. The error dynamics can be written as
\begin{equation}\label{E68}
\begin{aligned}
\dot{e}_1 &= -L_1 |e_1|^{\frac{2}{3}}\text{sign}(e_1)+e_2 \\
\dot{e}_2 &= -L_2|e_1|^{\frac{1}{3}}\text{sign}(e_1)+2180e_1 -\int_0^t L_3 \text{sign}(e_1) \mathrm{d}\tau +w(t,x)
\end{aligned}
\end{equation}
By means of transformation $e_{3}=-\int_0^t L_3 \text{sign}(e_1) \mathrm{d}\tau +w(t,x)$ and with the assumption 4, the system \eqref{E68}can be written as
\begin{equation}\label{E69}
\begin{aligned}
\dot{e}_1 &= -L_1 |e_1|^{\frac{2}{3}}\text{sign}(e_1)+e_2 \\
\dot{e}_2 &= -L_2|e_1|^{\frac{1}{3}}\text{sign}(e_1)+2180e_1+e_{3}\\
\dot{e}_3 &=-L_3 \text{sign}(e_1)+\dot{w}(t,x)
\end{aligned}
\end{equation}
In [46], it has been already proved that the \eqref{E69} is finite time stable and by selecting the appropriate gains $L_{1}$, $L_{2}$ and $L_{3}$ it can be guaranteed that $e_{1}$, $e_{2}$ and $e_{3}$ will converge to zero in finite time [47].\\
Now, the design of RTMF based on STA can be formulated with the estimated velocity $\hat{x}_2$. Using relation Using \eqref{E36}, \eqref{G}, \eqref{E50} and \eqref{E67} we have
\begin{align}
 \dot{z}_1 &= x_2-343000\,x_{r_2} \label{E70} \\
 \dot{z}_2 &=\dot{\hat{x}}_2-343000\dot{x}_{r_2} \nonumber\\
 &=L_2|e_1|^{\frac{1}{3}}\text{sign}(e_1) +\int_0^t L_3 \text{sign}(e_1) \mathrm{d}\tau \nonumber \\& +2180\,\hat{x}_1 -3518.85 \, u(t)-343000 x_{r_3} \label{E71}
 \end{align}
 Now choosing a sliding surface same as \eqref{E53}
 \begin{equation} \hat{s}=c_1z_1+ \hat{z}_2 \label{E72} \end{equation}
 Taking derivative of \eqref{E72} and using \eqref{E70} and \eqref{E71} we get
\begin{align}
 \dot{\hat{s}} &= c_1 \left(\hat{x}_2 +e_2 \right)-343000c_1  x_{r_2}+L_2|e_1|^{\frac{1}{3}}\text{sign}(e_1) \nonumber \\& +2180 \hat{x}_1-3518.85u(t) + \int_0^t L_3 \text{sign}(e_1) \mathrm{d}\tau -343000 x_{r_3} \label{E73}
 \end{align}
Now using \eqref{E57} and choosing $c_1=1$ we can write \eqref{E73} as
\begin{align}
\dot{\hat{s}} &= \left( \hat{x}_2 +e_2 \right)-343000 x_{r_2}-7.47 \times 10^8 \, x_{r_1} \nonumber \\& +8885 \, x_{r_3}+ L_2|e_1|^{\frac{1}{3}}\text{sign}(e_1) +2180 \hat{x}_1 \nonumber \\&-3518.85\,v(t) + \int_0^t L_3 \text{sign}(e_1) \mathrm{d}\tau \label{E74}
\end{align}
Now choosing $v(t)$ as
\begin{align}
v(t) &= -\frac{1}{3518.85} \bigg{(} - \hat{x}_2+343000 x_{r_2} + 7.47 \times 10^8 x_{r_1}\nonumber \\&  - 8885 x_{r_3} - L_2 |e_1|^{\frac{1}{3}} \text{sign}(e_1)- \int_0^t L_3 \text{sign}(e_1) \mathrm{d}\tau \nonumber \\& -2180 \hat{x}_1 -\lambda_1 |\hat{s}|^{\frac{1}{2}}\text{sign}(e_1) -\int_0^t \lambda_2 \text{sign}(\hat{s}) \mathrm{d}\tau \bigg{)} \label{E75}
\end{align}
Where $\lambda_1$ and $\lambda_2$ are controller design parameters, then we can write \eqref{E74} as
\begin{align}
\dot{\hat{s}} &=e_2 -\lambda_1 |\hat{s}|^{\frac{1}{2}}\text{sign}(e_1) -\int_0^t \lambda_2 \text{sign}(\hat{s}) \mathrm{d}\tau \label{E76}
\end{align}
Now using the fact $\dot{z}_1=\hat{z}_2$, \eqref{E72} and \eqref{E76} the closed loop system can be represented in $(z_1,\hat{s})$ coordinates as follows
\begin{align}
\dot{z}_1 &= -z_1+ \hat{s} \nonumber \\
\dot{\hat{s}} &=e_2 -\lambda_1 |\hat{s}|^{\frac{1}{2}}\text{sign}(e_1) -\int_0^t \lambda_2 \text{sign}(\hat{s}) \mathrm{d}\tau \label{E77}
\end{align}
So the overall closed loop dynamics with controller-observer together can be represented as
\begin{align}
\Phi_1 :& \begin{cases}
    \dot{z}_1=-{z}_1+\hat{s}  & \\
    \dot{\hat{s}}=e_2-\lambda_1|\hat{s}|^{1/2}\text{sign}(\hat{s})-\int_0^t \lambda_2 \text{sign}(\hat{s})\,\mathrm{d}\tau &
\end{cases} \nonumber \\
\Lambda_1 :& \begin{cases}
  \dot{e}_1 &= -L_1 |e_1|^{\frac{2}{3}}\text{sign}(e_1)+e_2 \\
\dot{e}_2 &= -L_2|e_1|^{\frac{1}{3}}\text{sign}(e_1)+2180e_1+e_{3}\\
\dot{e}_3 &=-L_3 \text{sign}(e_1)+\dot{w}(t,x)
\end{cases} \label{E77}
\end{align}
It has been already discussed that the observer error for the system $\Lambda_1$ goes to zero in finite time. In [45] it has been established that, the trajectories of the system $\Phi_1$ cannot escape in finite time. By substituting $e_{2}=0$, the closed-loop system can be written as
\begin{align}\label{E62}
\begin{split}
    \dot{z}_1&=-z_1+\hat{s}  \\
    \dot{\hat{s}}&=-\lambda_1|\hat{s}|^{1/2}\text{sign}(\hat{s})+\xi \\
    \dot{\xi}&=-\lambda_2\text{sign}(\hat{s})
\end{split}
\end{align}
In \eqref{E62}, the lower two equations are of STA and by selecting appropriate gains $\lambda_1>0$ and $\lambda_2>0$ it can be ensured that $\hat{s}=\dot{\hat{s}}=0$ in finite time which helps to represt the reduced order system as follows
\begin{align}
\dot{z}_1 &= -z_1 \label{E78} \\
\hat{z}_2 &= - z_1 \label{E79}
\end{align}
Therefore $z_1$ and ${z}_2$ goes to zero asymptotically which in turn guarantees the asymptotic convergence of output tracking error $e(t)$.
\section{Simulation and experimental results}
This section provides the experimental results for MagLev system for two cases a) RTMF controller based on STA with STO and b) RTMF controller based on STA with HOSMO. The simulation and experiment with MagLev system has been performed for these two cases. The external disturbance voltage signal $w(t)=5\text{sin}(t)$ is injected externally to perturb the plant. The output tracking performance of MagLev system is presented for sinusoidal and trapezoidal input for both cases which is mentioned above. The initial condition for the model is taken as $x_{r}(0) = [1 \times 10^{-5}~~0~~0]^{\top}$.

\subsection{RTMF controller based on STA with STO}
\begin{figure}[h!]
   \centering
 \includegraphics[width=3.5in]{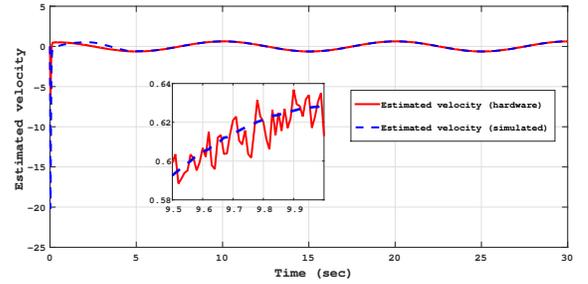}
 \caption{Simulated and experimental estimation of velocity using STO (sinusoidal)}
 \label{fig2}
\end{figure}

\begin{figure}[h!]
   \centering
 \includegraphics[width=3.5in]{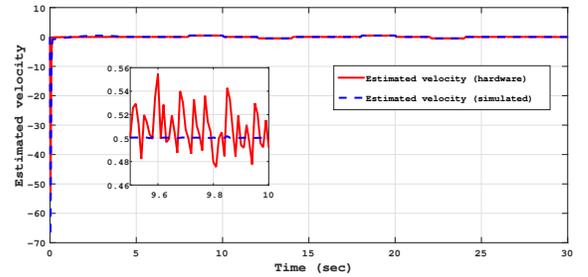}
 \caption{Simulated and experimental estimation of velocity using STO (trapezoidal)}
 \label{fig3}
\end{figure}

The control strategy based on STO developed in Section 5.2 is implemented for MagLev system. The STO gains are chosen as $K_1=50$ and $K_2=400$ and gains of RTMF controller based on STA  are chosen as $\lambda_1=\lambda_2=10$. Fig.~\ref{fig2} and Fig.~\ref{fig3} shows the estimated velocity for the inputs sinusoidal and trapezoidal respectively with STO.

\begin{figure}[h!]
   \centering
 \includegraphics[width=3.5in]{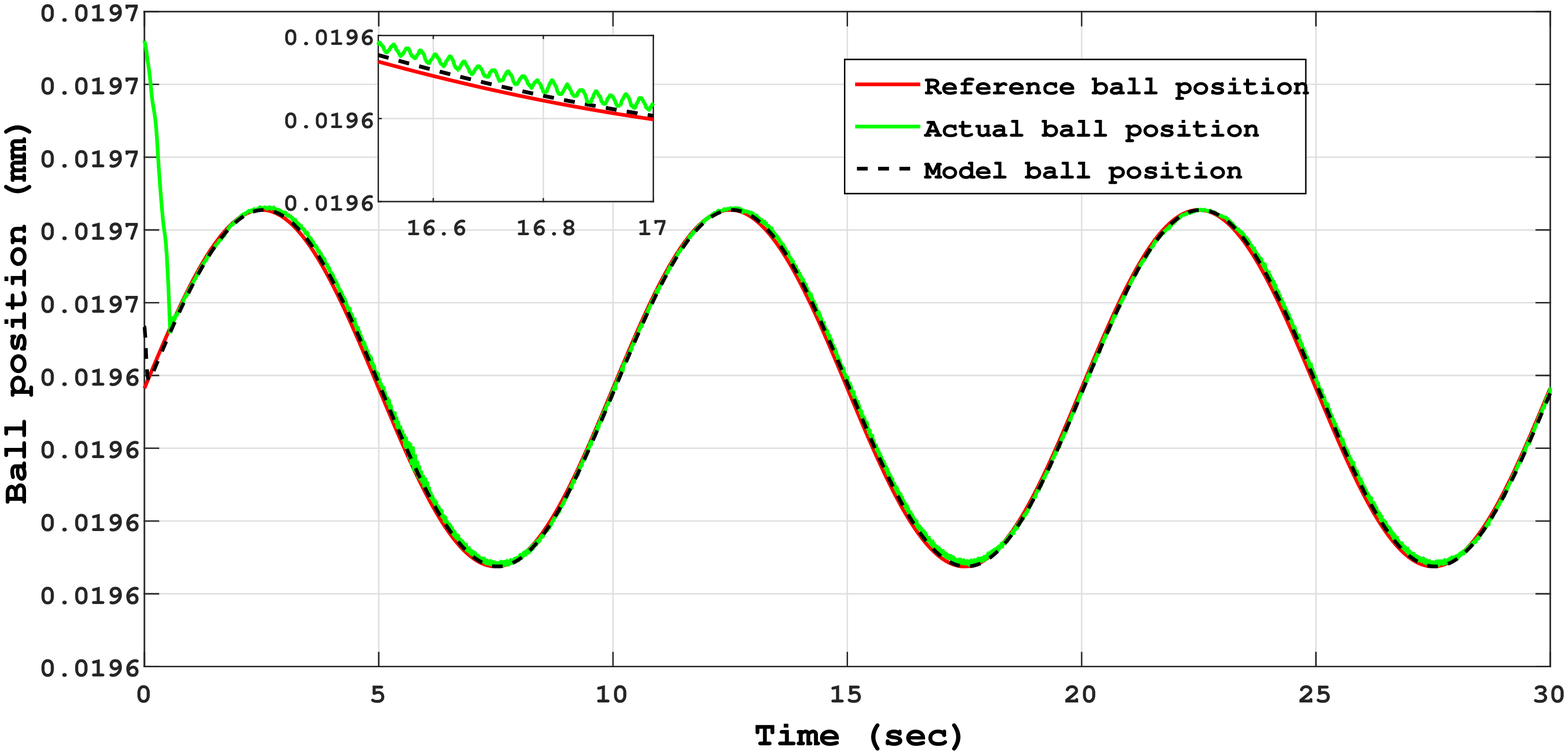}
 \caption{Tracking performance with sinusoidal input (experimental)}
 \label{fig4}
\end{figure}

\begin{figure}[h!]
   \centering
 \includegraphics[width=3.5in]{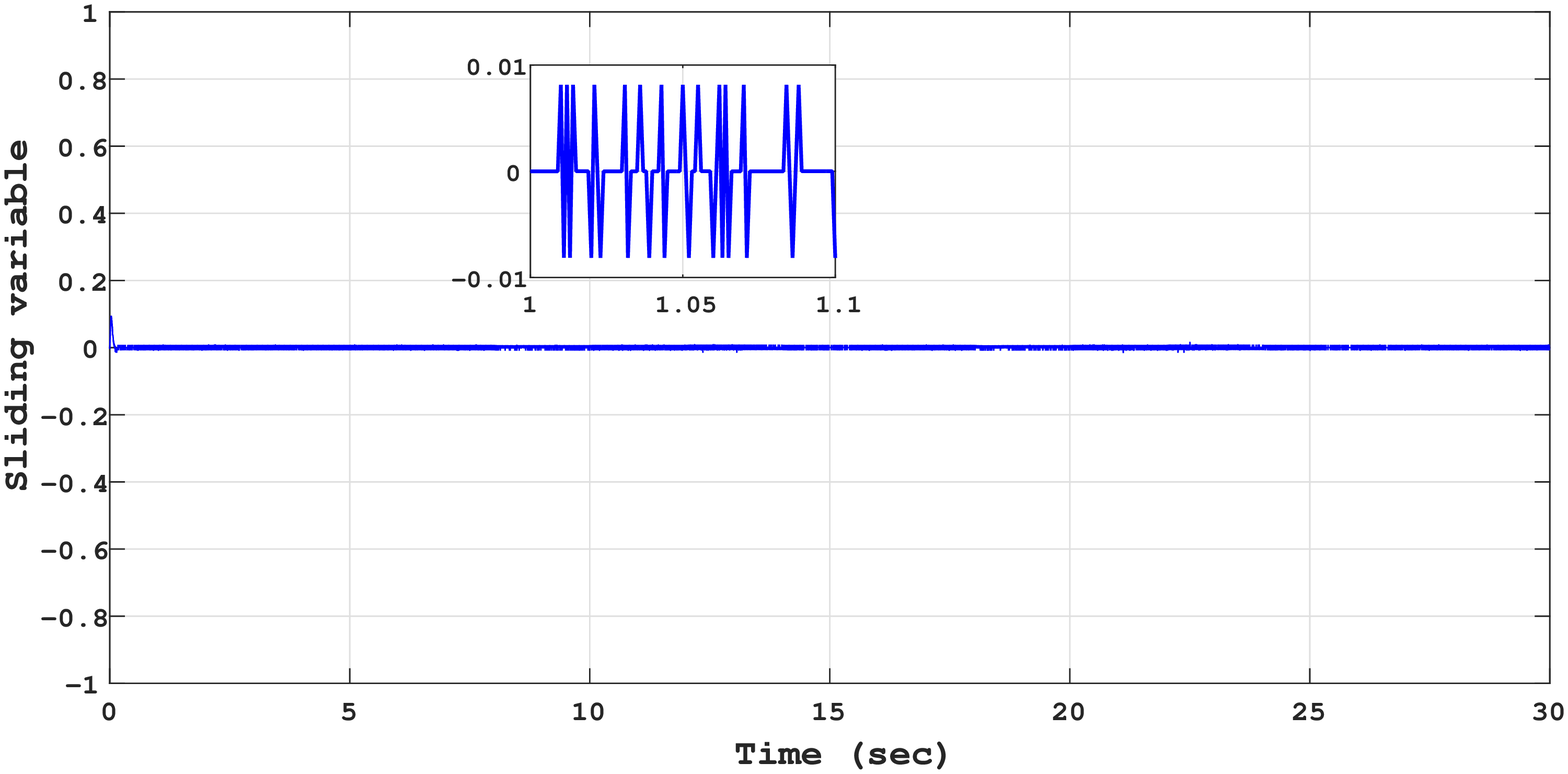}
 \caption{Sliding surface (experimental)}
 \label{fig5}
\end{figure}

\begin{figure}[h!]
   \centering
 \includegraphics[width=3.5in]{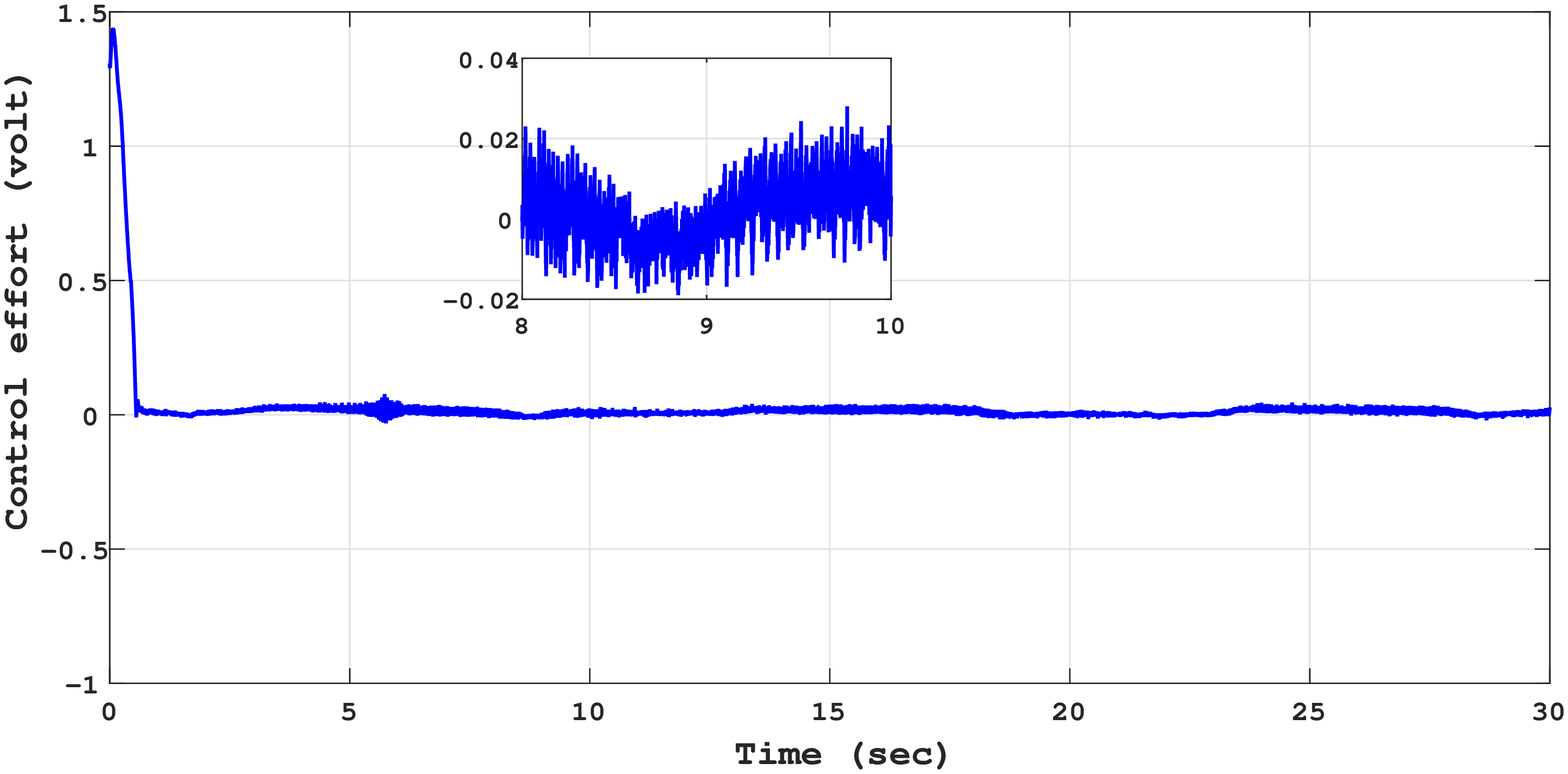}
 \caption{Control effort (experimental)}
 \label{fig6}
\end{figure}

The experimental results for sinusoidal tracking performance, sliding surface and control input is presented in Fig.~\ref{fig4}, Fig.~\ref{fig5} and Fig.~\ref{fig6} respectively.

\begin{figure}[h!]
   \centering
 \includegraphics[width=3.5in]{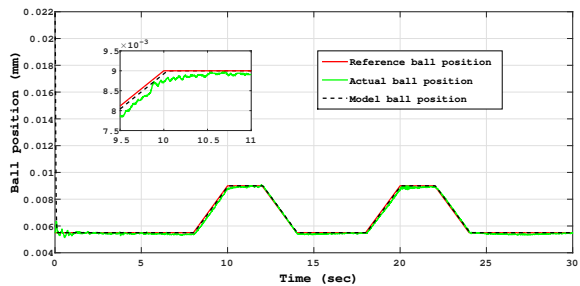}
 \caption{Tracking performance with trapezoidal input (experimental)}
 \label{fig7}
\end{figure}

\begin{figure}[h!]
   \centering
 \includegraphics[width=3.5in]{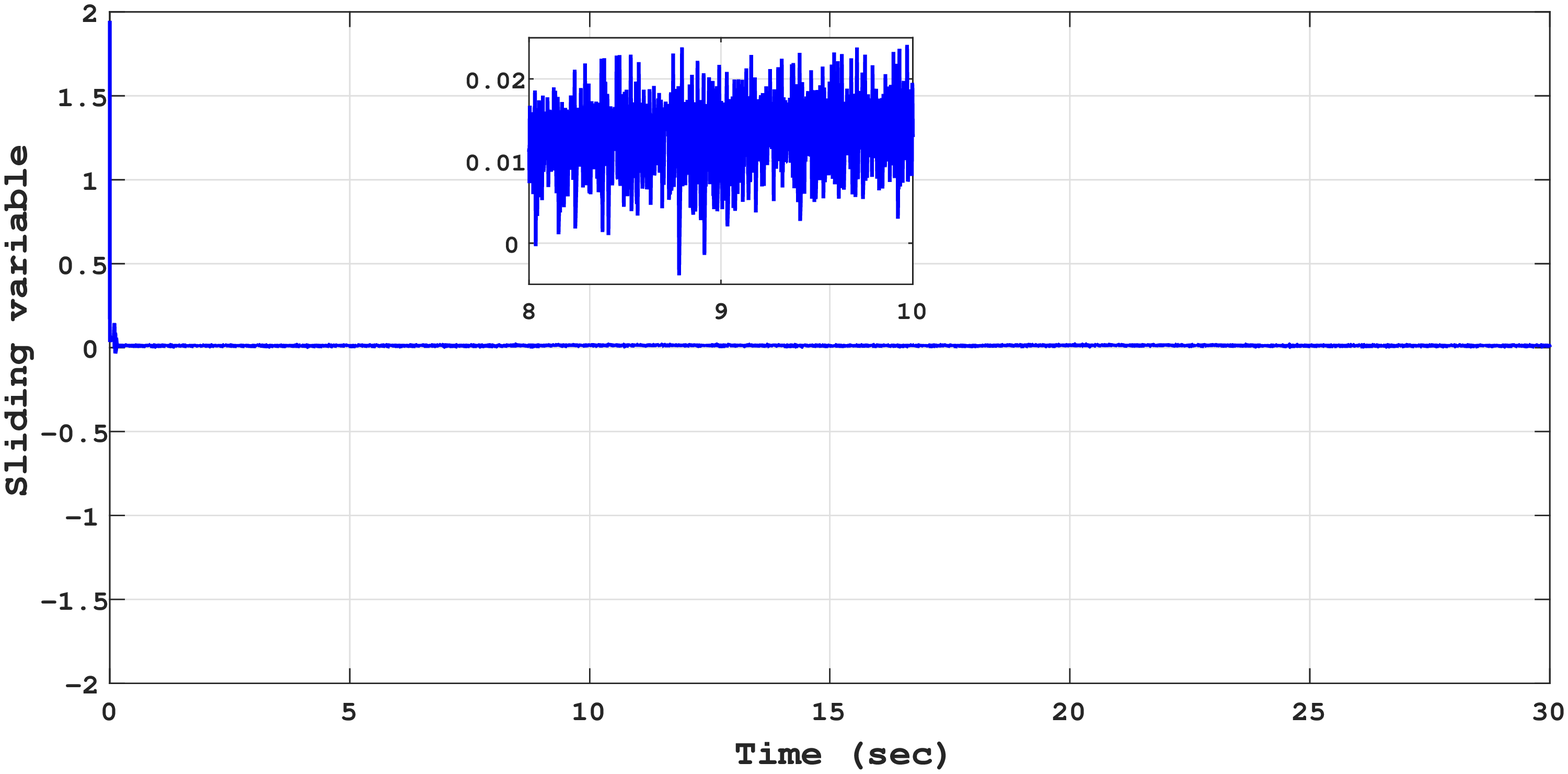}
 \caption{Sliding surface (experimental)}
 \label{fig8}
\end{figure}

\begin{figure}[h!]
   \centering
 \includegraphics[width=3.5in]{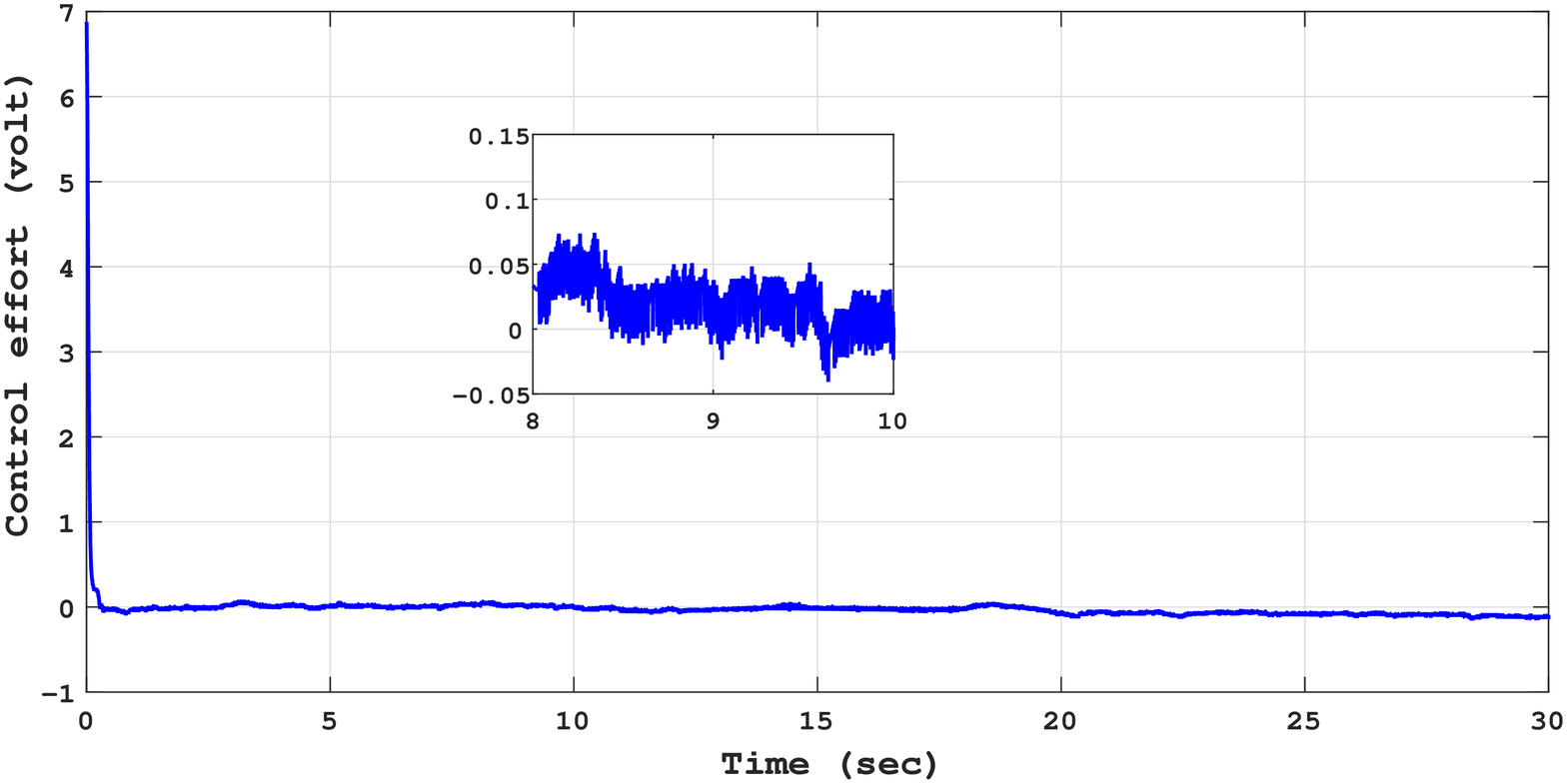}
 \caption{Control effort (experimental)}
 \label{fig9}
\end{figure}

Similarly the experimental results for trapezoidal tracking performance, sliding surface and control input is presented in Fig.~\ref{fig7}, Fig.~\ref{fig8} and Fig.~\ref{fig9} respectively.

\subsection{RTMF controller based on STA with HOSMO}
\begin{figure}[h!]
   \centering
 \includegraphics[width=3.5in]{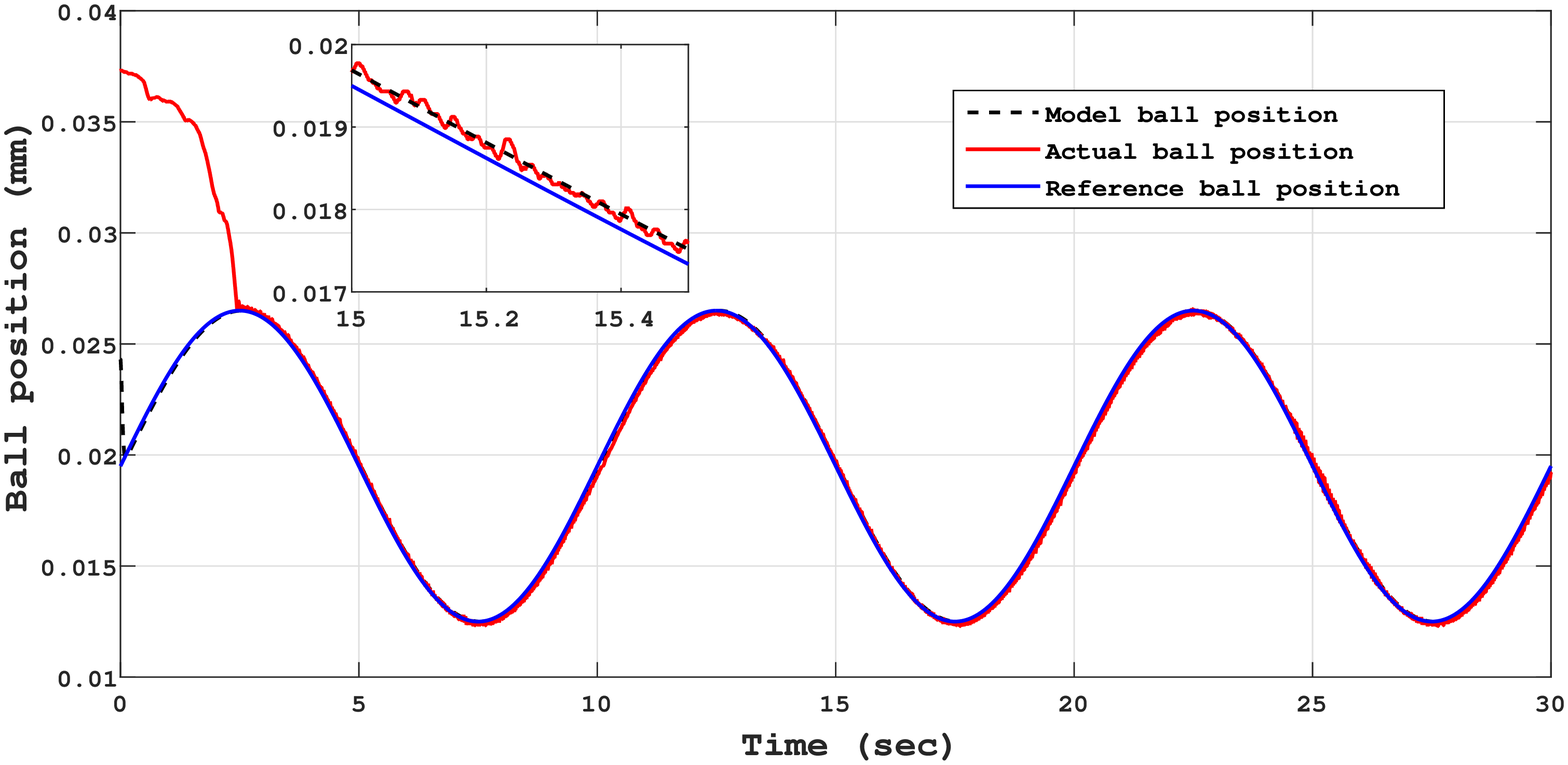}
 \caption{Tracking performance with sinusoidal input (experimental)}
 \label{fig10}
\end{figure}

\begin{figure}[h!]
   \centering
 \includegraphics[width=3.5in]{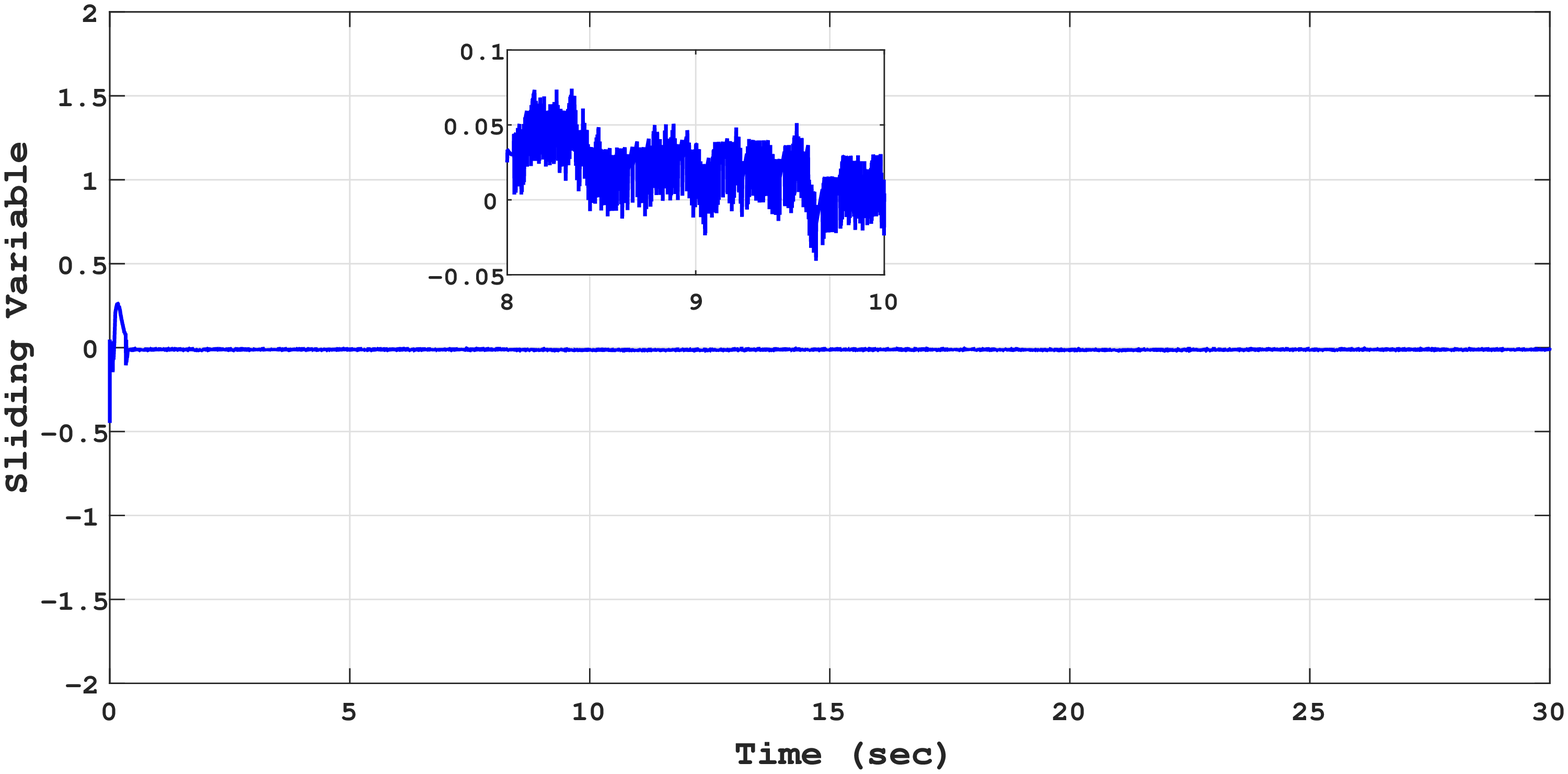}
 \caption{Sliding surface (experimental)}
 \label{fig11}
\end{figure}

\begin{figure}[h!]
   \centering
 \includegraphics[width=3.5in]{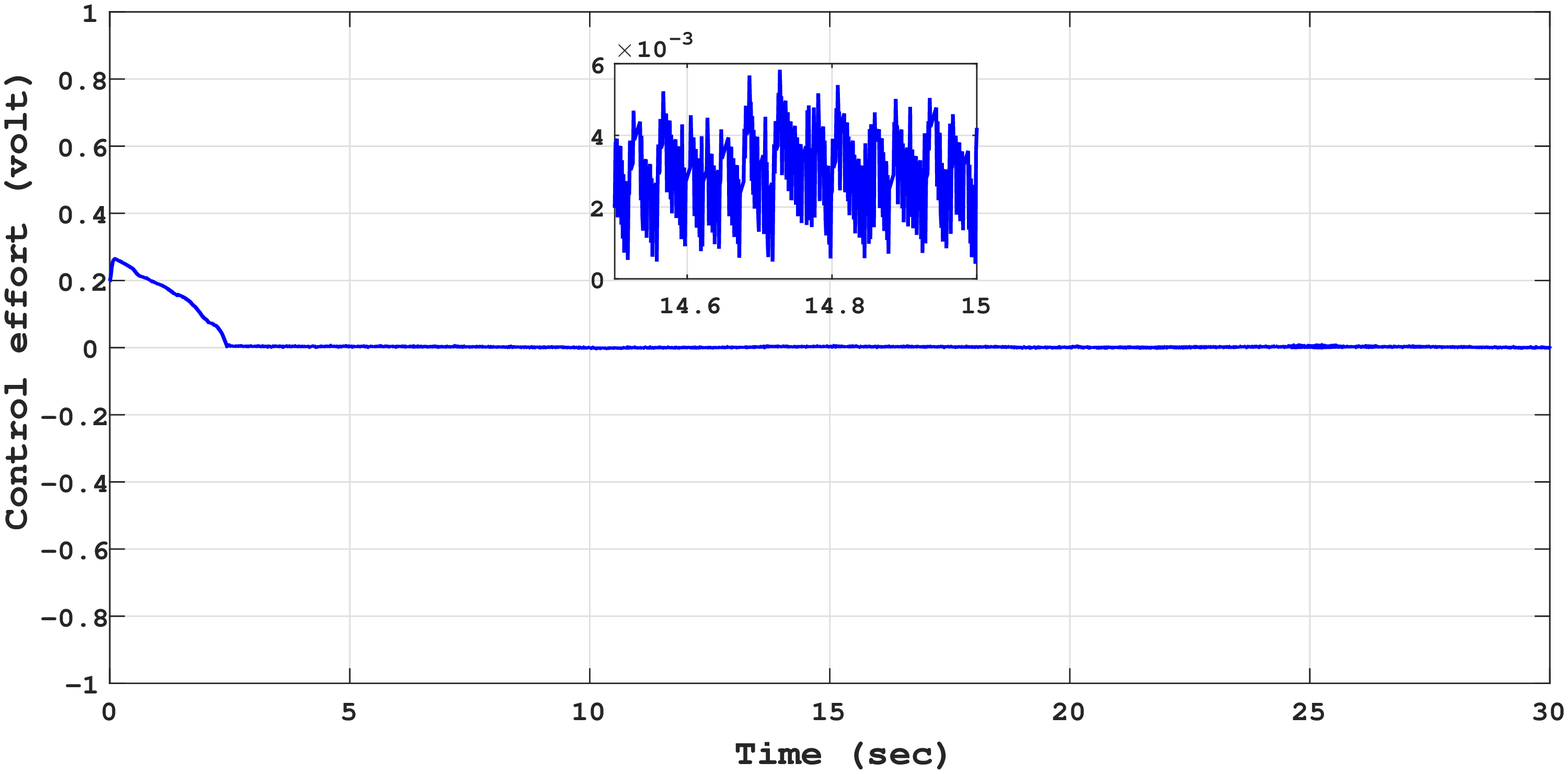}
 \caption{Control effort (experimental)}
 \label{fig12}
\end{figure}

This section is devoted to present the experimental results for the RTMF controller based on STA with HOSMO. To estimate the velocity of the ball of the MagLev system HOSMO is implemented with gains : $L_1=35, L_2=100$ and $L_3=600$. The gains of the RTMF controller based on STA is chosen as $\lambda_1=\lambda_2=15$. The estimation performance of HOSMO is found similar to STO and the results are not shown because of the space constraint. The experimental results of this controller-observer combination for sinusoidal tracking, sliding surface and control input is presented in Fig.~\ref{fig10}, Fig.~\ref{fig11} and Fig.~\ref{fig12} respectively.

\begin{figure}[h!]
   \centering
 \includegraphics[width=3.5in]{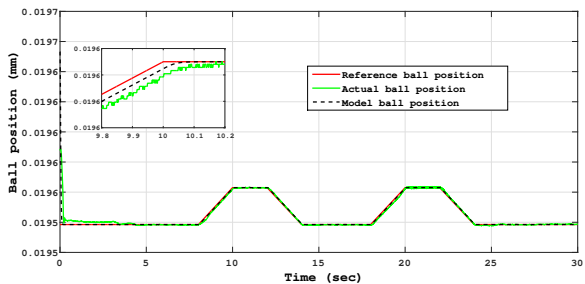}
 \caption{Tracking performance with trapezoidal input (experimental)}
 \label{fig13}
\end{figure}

\begin{figure}[h!]
   \centering
 \includegraphics[width=3.5in]{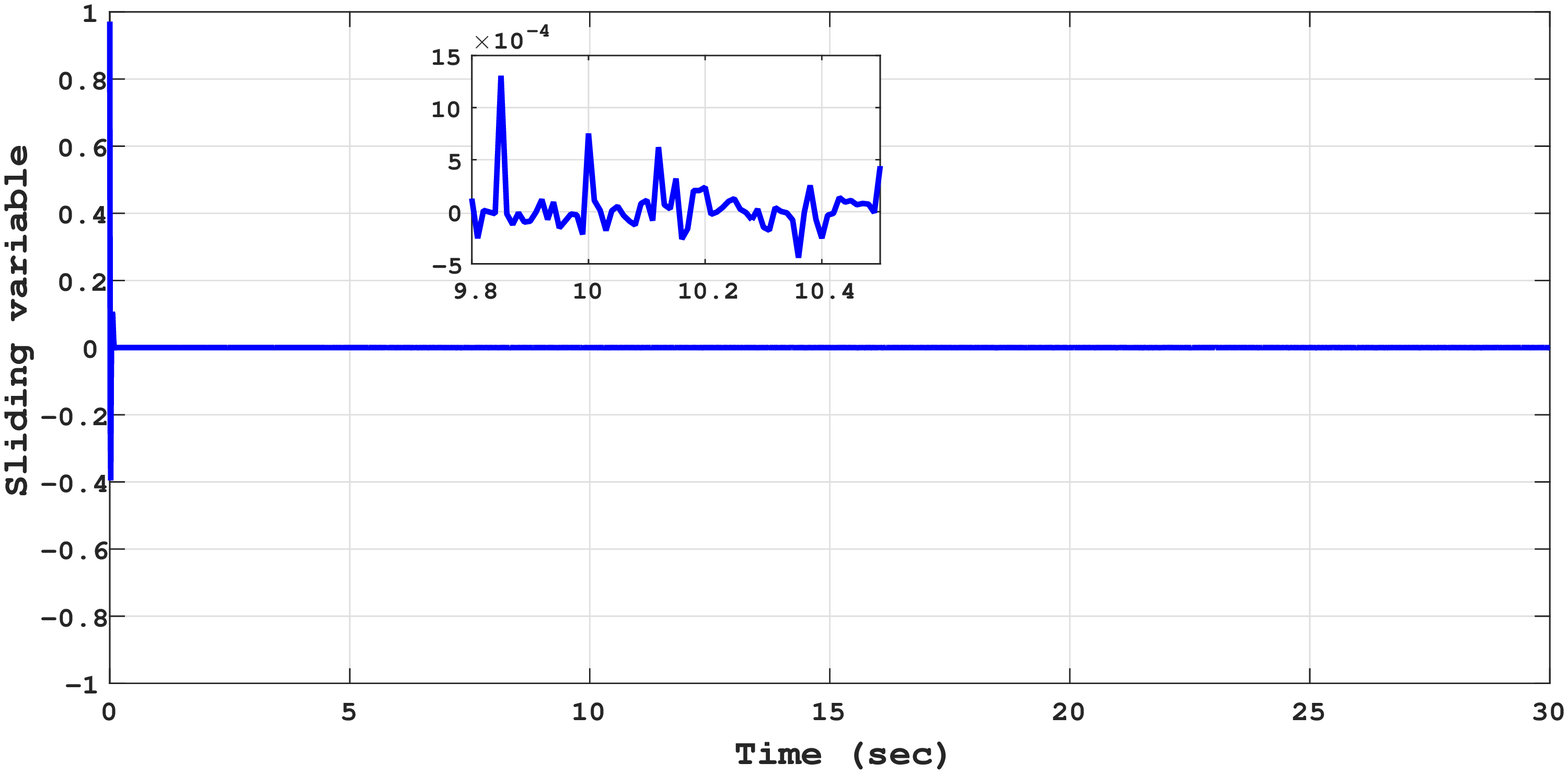}
 \caption{Sliding surface (experimental)}
 \label{fig14}
\end{figure}

\begin{figure}[h!]
   \centering
 \includegraphics[width=3.5in]{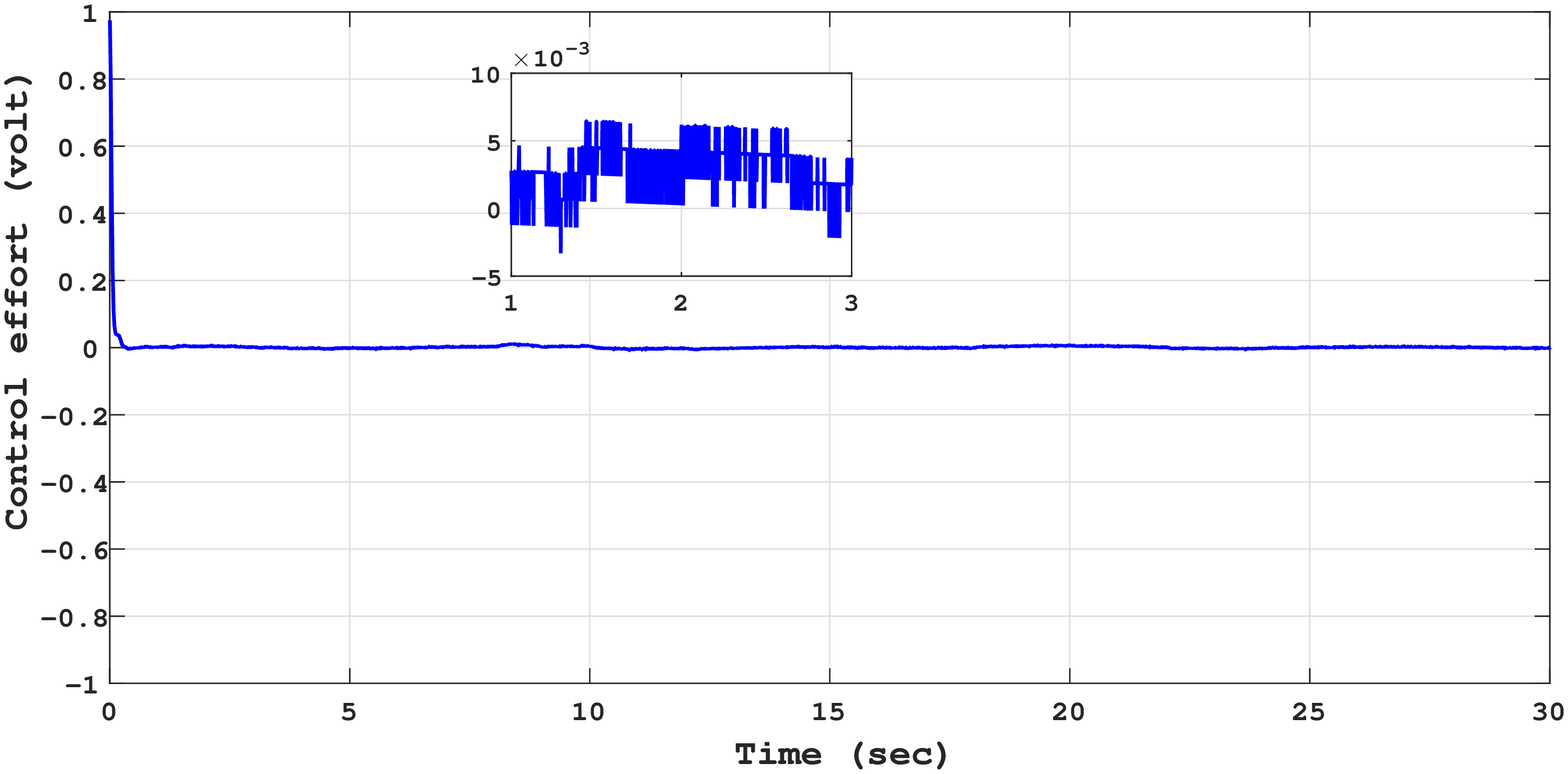}
 \caption{Control effort (experimental)}
 \label{fig15}
\end{figure}

In a similar manner the  experimental results for trapezoidal tracking performance, sliding surface and control input is presented in Fig.~\ref{fig13}, Fig.~\ref{fig14} and Fig.~\ref{fig15} respectively.

\par Comparing the control plots Fig.~\ref{fig6}, Fig.~\ref{fig9} and Fig.~\ref{fig12}, Fig.~\ref{fig15}, it is easy to conclude that the control signal is more smoother in case of RTMF controller based on STA with HOSMO compared to RTMF controller based on STA with STO. It is apparent form the sliding surface plots (in zoom plot), the improved precision is achieved in case of RTMF controller based on STA with HOSMO compared to RTMF controller based on STA with STO.

\section{Conclusion}\label{sec10}
In this paper, RTMF controller based on STA is proposed for uncertain LTI systems. To validate the performance of the proposed controller, it is designed and implemented in real time MagLev system in the presence of disturbance. The implementation of this controller for MagLev requires knowledge of both the system states. However only ball position is available for measurement. To avoid the implementation difficulties, the RTMF controller based on STA
is designed and implemented with STO. But it has been realized that the continuous control for the chosen sliding surface is not possible to achieve with  RTMF controller based on STA with STO. To overcome this issue, a RTMF controller based on STA with HOSMO is designed and implemented for MagLev system. It is clear from the obtained results that the control scheme based on STA provides excellent tracking of time varying signals for both cases during which the system is completely insensitive to disturbances acting on it. It has been also observed that smoother control action with improved precision of sliding variable is achieved in case of RTMF controller based on STA with HOSMO.

\end{document}